\documentclass[11pt,reqno]{amsart}

\usepackage{amsmath,amssymb}

\def\demo{\noindent{\bf Proof. }}
\def\sqr#1#2{{\vcenter{\hrule height.#2pt
            \hbox{\vrule width.#2pt height#1pt \kern#1pt
                    \vrule width.#2pt}
            \hrule height.#2pt}}}
\def\square{\mathchoice\sqr64\sqr64\sqr{4}3\sqr{3}3}
\def\QED{\hfill$\square$}

\def\lar{\longrightarrow}
\def\rar{\rightarrow}

\def\Hom{{\rm Hom}}
\def\gr{{\rm gr}}

\newtheorem{Theorem}{Theorem}[section]

\newtheorem{Corollary}[Theorem]{Corollary}
\newtheorem{Proposition}[Theorem]{Proposition}

\newtheorem{Remark}[Theorem]{Remark}
\newtheorem{Example}[Theorem]{Example}

\newtheorem{Question}[Theorem]{Question}

\begin{document}

\baselineskip=15pt

\title[Multiplicity of the special fiber of blowups]{\bf
Multiplicity of the special fiber of blowups}

\author[A. Corso, C. Polini and W.V. Vasconcelos]
{Alberto Corso \and Claudia Polini \and Wolmer V. Vasconcelos}

\address{Department of Mathematics, University of Kentucky,
Lexington, Kentucky 40506 - USA} \email{corso@ms.uky.edu}

\address{Department of Mathematics, University of Notre Dame,
Notre Dame, Indiana 46556 - USA} \email{cpolini@nd.edu}

\address{Department of Mathematics, Rutgers University, New Brunswick,
New Jersey 08901 - USA} \email{vasconce@math.rutgers.edu}

\thanks{
AMS 2000 {\em Mathematics Subject Classification}. Primary 13A30,
13B21, 13D40; Secondary 13H10, 13H15. \newline\indent Part of this
work was done while the authors were members of the Mathematical
Sciences Research Institute (Berkeley) during the Fall 2002. They
would like to thank that institution as well as the organizers of
the year long program in Commutative Algebra for the stimulating
atmosphere that this event provided. The last two authors also
gratefully acknowledge partial support from the NSF}

\begin{abstract}
Let $(R, {\mathfrak m})$ be a Noetherian local ring and let $I$ be
an ${\mathfrak m}$-primary ideal. In this paper we give sharp
bounds on the multiplicity of the special fiber ring ${\mathcal
F}$ of $I$ in terms of other well-known invariants of $I$. A
special attention is then paid in studying when equality holds in
these bounds, with a particular interest in the unmixedness or, better,
the Cohen-Macaulayness of ${\mathcal F}$.
\end{abstract}

\vspace{-0.1in}

\maketitle

\vspace{-0.2in}

\section{Introduction}

Let $(R, {\mathfrak m})$ be a Noetherian local ring with dimension
$d>0$ and infinite residue field and let $I$ be an $R$-ideal.
Three algebras, collectively referred to as {\it blowup algebras of $I$},
play an important role in the process of blowing up the variety
${\rm Spec}(R)$ along the subvariety $V(I)$: namely, the {\it Rees
algebra} ${\mathcal R}(={\mathcal R}(I)=R[It])$, the {\it associated
graded ring} ${\mathcal G}(={\mathcal G}(I))$, and the {\it
special fiber ring} ${\mathcal F}(={\mathcal F}(I))$ of $I$
\[
{\mathcal R}= \bigoplus_{m=0}^{\infty} I^mt^m, \qquad {\mathcal G} =
{\mathcal R}/I{\mathcal R}, \qquad {\mathcal F} ={\mathcal
R}/{\mathfrak m}{\mathcal R}.
\]

The blowup algebras are also extensively used as the means to
examine diverse properties of the ideal $I$. Therefore, much
attention has been paid in the past to find under which
circumstances these objects have a good structure. In this article
we focus on the special fiber ring ${\mathcal F}$, as there is
lack of knowledge about its properties in comparison to the other two
blowup algebras of $I$. From an algebraic point of view,
${\mathcal F}$ yields $($asymptotic$)$ information about the ideal
$I$. For example, its Hilbert function is the numerical function
that measures the growth of the minimal number of generators
$\mu(I^m)$ of the powers of $I$. For $m \gg 0$ this function is a
polynomial in $m$ of degree ${\rm dim}\, {\mathcal F}-1$, whose
leading coefficient $f_0=f_0(I)$ is called the {\it multiplicity} $($or
{\it degree}$)$ of the special fiber ring ${\mathcal F}$.
Another important data attached to ${\mathcal F}$ is its
$($Krull$)$ dimension, dubbed the {\it analytic spread} $\ell$ of
$I$. It is bounded below by the height $g$ of $I$ and bounded
above by the dimension of $R$. It coincides with the minimal
number of generators of any minimal reduction of $I$. A
minimal reduction --- a notion that has been crucial in the study
of the Rees algebra of an ideal, as it carries most of the
information about the original ideal but, in general, with fewer
generators --- arise from a Noether normalization of ${\mathcal
F}$.

From a more geometric perspective, ${\rm Proj}({\mathcal F})$
corresponds to the fiber over the closed point of the blowup of
Spec($R$) along $V(I)$. When $R$ is a standard graded domain over
a field $k$ and $I$ is the $R$-ideal generated by forms $a_0,
\ldots, a_n$ of the same degree, then ${\mathcal F}$ describes the
homogeneous coordinate ring of the image of the rational map ${\rm
Proj}(R) \dashrightarrow {\mathbb P}_k^n$ given by $(a_0, \ldots,
a_n)$. As a special case this construction yields homogeneous
coordinate rings of Gauss images and of secant varieties. More
recently, special fiber rings also find an application in the
theory of evolutions, due to Mazur \cite{Mazur} and inspired by
the work of Wiles \cite{W} on semistable curves. Indeed, in this
context, H\"ubl and Huneke \cite{HuHu} were the first to use the
special fiber ring ${\mathcal F}$ in studying the issue of
the integral closedness of the ideal ${\mathfrak m}I$. See
\cite{CGPU} for subsequent work.

We now describe the content of the paper. In Section 2 we prove
our main results, which deal with an ${\mathfrak m}$-primary ideal
$I$ of a local Cohen-Macaulay ring $R$. In Theorem~\ref{sally} we
show that the multiplicity $f_0$ of the special fiber ring
${\mathcal F}$ satisfies the inequality
\[
f_0 \leq  e_1 - e_0 + \lambda(R/I) + \mu(I) - d + 1,
\]
where $e_0=e_0(I)$ and $e_1=e_1(I)$ are the first two coefficients of the
Hilbert polynomial of $I$ and $\lambda(\cdot)$ denotes the
length function. We recall that for $m \gg 0$ the Hilbert-Samuel
function of an ${\mathfrak m}$-primary ideal $I$ of a local
Cohen-Macaulay ring $R$ --- that is the numerical function that
measures the growth of $\lambda(R/I^m)$ for all $m \geq 1$ --- is
a polynomial in $m$ of degree $d$, whose coefficients $e_0, e_1,
\ldots, e_d$ are called the Hilbert coefficients of $I$. The
leading term $e_0$ of this polynomial is called the multiplicity
$($of the associated graded ring$)$ of $I$: It is the only
coefficient which is geometrically well understood as it equals
$\lambda(R/J)$, for any minimal reduction $J$ of $I$.

Even though equality in the above bound does not assure the
Cohen-Macaulayness of ${\mathcal F}$ $($see
Example~\ref{minors}$)$, nevertheless it yields useful information
about the associated primes of ${\mathcal F}$. In fact, in
Theorem~\ref{cmfif} we show that if $f_0 = e_1-e_0 + \lambda(R/I)+
\mu(I)-d+1$, then ${\mathcal F}$ is unmixed.
In particular the ideals ${\mathfrak m}I^m$ are integrally closed
for all $m$, whenever ${\mathcal R}$ is normal $($see
Corollary~\ref{integrality}$)$. An immediate consequence of
Theorem~\ref{cmfif} is that if ${\rm depth}\, {\mathcal G} \geq d-1$
then ${\mathcal F}$ is Cohen-Macaulay $($see Corollary~\ref{depthd-1}$)$.

Another motivation in finding bounds on $f_0$ in terms of natural data
attached to the ideal $I$ is based on the general philosophy that
a good control on the multiplicity $f_0$ coupled with depth
information on ${\mathcal F}$ allow to bound the reduction number $r$ of
the ideal. We recall that $J$ is a {\it reduction} of $I$ if
$I^{r+1}=JI^r$ for some non-negative integer $r$ or, equivalently, if
the inclusion of Rees algebras ${\mathcal R}(J) \hookrightarrow
{\mathcal R}(I)$ is module finite \cite{NR}. The least such $r$
is called the reduction number $r_J(I)$ of $I$ with respect to $J$.
One then defines the {\it reduction number} of $I$ to be the least
$r_J(I)$, where $J$ varies over all minimal reductions of $I$.
A reduction is said to be {\it minimal} if it is minimal with respect
to containment. In other words, the reduction number of an ideal is
a key control element of the blowup algebras, often measuring the interplay
among the other invariants of the ideal.

For example, if $I$ is generated by forms of the same degree then
${\mathcal F}$ is a domain and the reduction number
$r$ of $I$ is bounded by $f_0$. Thus, in Section 3 we use the
bounds on $f_0$ established earlier in the paper to obtain bounds on the
reduction number of an ${\mathfrak m}$-primary ideal $I$. For
instance, if the residue field $R/{\mathfrak m}$ has
characteristic $0$ and $f_0 = e_1 - e_0 + \lambda(R/I)+
\mu(I)-d+1$ we conclude that $r \leq e_1 - e_0 +\lambda(R/I) +
\mu(I) - d$ $($see Corollary~\ref{firstcor}$)$. If, in addition,
${\mathcal F}$ is Cohen-Macaulay this bound can be improved and in
Corollary~\ref{coro-rossi} we show that $r \leq e_1 - e_0
+\lambda(R/I) + 1$. This latter bound has been shown without any
additional assumption by Rossi \cite{R} if the dimension of the
ring is at most 2. It seems to hold in full generality.

In Section 4 we strengthen the bound on $f_0$ obtained in Section
1. With the same assumption as in Theorem~\ref{sally}, we first
prove that the multiplicity $f_0$ of the special fiber ring
${\mathcal F}$ of $I$ satisfies the tighter bound $f_0 \leq e_1 -
e_0 + \lambda(R/\widetilde{I}) + \mu(\widetilde{I}) - d +1$, where
$\widetilde{I}$ denotes the Ratliff-Rush closure of $I$ $($see
Corollary~\ref{rr-closure}$)$. We recall that if $I$ is an
${\mathfrak m}$-primary ideal containing a regular element, then
the {\it Ratliff-Rush closure} $\widetilde{I}$ of an ideal $I$ is
the largest ideal containing $I$ with the same Hilbert polynomial
as $I$ \cite{RR}. We can actually improve the estimate on $f_0$ by
giving a different derivation, that we carry out in
Theorem~\ref{s2ofideal}, of the construction, originally due to
Shah \cite{Shah}, of a canonical sequence of ideals containing $I$
with partially identical Hilbert polynomials to the one of $I$. In
Corollary~\ref{e1-closure} we then show that $f_0 \leq e_1 - e_0 +
\lambda(R/\check{I}) + \mu(\check{I}) - d + 1$, where $\check{I}$
is the degree one component of the $S_2$-ification of the Rees
algebra of $I$. The ideal $\check{I}$ can also be characterized as
the largest ideal containing $I$ with the same $e_0$ and $e_1$ as
$I$ \cite{Ciu}. The crucial issue is then to show that the
multiplicity $f_0$ of the special fiber ring of $I$ is unchanged
when passing to the special fiber ring of $\check{I}$.

The main idea, behind the proofs of Theorem~\ref{sally} and
Theorem~\ref{cmfif} is the use of the Sally module $S_J(I)$ of $I$
with respect to a minimal reduction $J$, which is defined by the
following exact sequence of ${\mathcal R}(J)$-modules
\[
0 \rightarrow I{\mathcal R}(J) \longrightarrow I{\mathcal R}(I)
\longrightarrow S_J(I) = \bigoplus_{m \geq 2}I^m/IJ^{m-1}
\rightarrow 0,
\]
introduced by Vasconcelos \cite{Vas}. In this paper, though, we
make use of a new approach to $S_J(I)$ by means of the following
exact sequence
\[
{\mathcal R}(J) \oplus  {\mathcal R}(J)^{n-\ell}[-1]
\longrightarrow {\mathcal R}(I) \longrightarrow S_J(I)[-1]
\rightarrow 0,
\]
where $n=\mu(I)$ and $\ell$ is the analytic spread of $I$. The
advantage of this approach is that it is suitable to study the
multiplicity of the special fiber ring ${\mathcal F}$ of any
ideal, not necessarily ${\mathfrak m}$-primary.
In addition, this approach avoids the customary technique in the
theory of Hilbert function of going modulo a superficial sequence.


\section{General bounds}

We start describing some general bounds on the multiplicity $f_0$
of the special fiber ring ${\mathcal F}$ of $I$. A special
attention is then paid in studying when equality holds in these
bounds, with a particular interest in the unmixedness or, better,
the Cohen-Macaulayness of ${\mathcal F}$.

\subsection{Upper bounds}
In \cite[2.4]{mrn} Vasconcelos has shown that the multiplicity
$f_0$ of the special fiber ring ${\mathcal F}$ of an ${\mathfrak
m}$-primary ideal $I$ of a local Cohen-Macaulay ring $R$ satisfies
the inequality
\begin{equation}\label{e0&e1}
f_0 \leq \min\{ e_0, \, e_1+1 \},
\end{equation}
where $e_0$ and $e_1$ are the first two coefficients of the Hilbert
polynomial of $I$. In the following theorem  we give a better bound on
$f_0$ of a different nature, using the structure of the Sally module
$S_J(I)$ of $I$ with respect to a minimal reduction $J$.

\begin{Theorem}\label{sally}
Let $(R, {\mathfrak m})$ be a local Cohen-Macaulay ring of
dimension $d>0$ and infinite residue field. Let $I$ be an
${\mathfrak m}$-primary ideal. Then the multiplicity $f_0$ of the
special fiber ring ${\mathcal F}$ of $I$ satisfies
\[
f_0 \leq e_1 - e_0 + \lambda(R/I) + \mu(I) - d + 1.
\]
\end{Theorem}
\demo Let $J$ be a minimal reduction of $I$ and write $I=(J, a_1,
\ldots, a_{n-d})$, where $n=\mu(I)$ and $d=\mu(J)=\ell$. We now
consider the Sally module $S_J(I)$ of $I$ with respect to $J$
\cite[2.1]{Vas}. However, we find it more suitable for our
purposes to approach $S_J(I)$ by means of a new exact sequence,
namely
\[
{\mathcal R}(J) \oplus  {\mathcal R}(J)^{n-d}[-1]
\stackrel{\varphi}{\longrightarrow} {\mathcal R}(I)
\longrightarrow S_J(I)[-1] \rightarrow 0,
\]
where $\varphi$ is the map defined by $\varphi(r_0, r_1, \ldots,
r_{n-d}) = r_0+r_1a_1t+\ldots+r_{n-d}a_{n-d}t$, for any element
$(r_0, r_1, \ldots, r_{n-d})\in {\mathcal R}(J) \oplus {\mathcal
R}(J)^{n-d}[-1]$. Tensoring the above exact sequence with
$R/{\mathfrak m}$ yields
\begin{equation}\label{tensor}
{\mathcal F}(J) \oplus {\mathcal F}(J)^{n-d}[-1] \longrightarrow
{\mathcal F}(I) \longrightarrow S_J(I)[-1] \otimes R/{\mathfrak m}
\rightarrow 0.
\end{equation}
As the three modules in $($\ref{tensor}$)$ have the same
dimension, we obtain the multiplicity estimate
\[
f_0 \leq \deg(S_J(I)[-1] \otimes
R/{\mathfrak m}) + \deg ({\mathcal F}(J) \oplus {\mathcal
F}(J)^{n-d}[-1]).
\]
Since $S_J(I) \otimes R/{\mathfrak m}$ is a homomorphic image of
$S_J(I)$, its multiplicity is bounded by the one of $S_J(I)$,
which --- according to \cite[3.3]{Vas} --- is
$e_1-e_0+\lambda(R/I)$. On the other hand, ${\mathcal F}(J) \oplus
{\mathcal F}(J)^{n-d}[-1]$ is a free ${\mathcal F}(J)$-module of
rank $n-d+1$. Thus, its multiplicity is $n-d+1$, since ${\mathcal
F}(J)$ is isomorphic to a ring of polynomials. \QED

\medskip

\begin{Proposition}
Let $(R, {\mathfrak m})$ be a local Cohen-Macaulay ring of
dimension $d>0$ and infinite residue field. Let $I$ be an
${\mathfrak m}$-primary ideal. We have that
\[
f_0 \leq e_1 - e_0 + \lambda(R/I) + \mu(I) - d + 1 \leq e_1 + 1.
\]
In particular, if $f_0=e_1+1$ then $I$ has {\rm minimal
multiplicity} in the sense of\/ {\rm S. Goto}, that is ${\mathfrak
m}I={\mathfrak m}J$ for any minimal reduction $J$ of $I$. If, in
addition, $R$ is a local Gorenstein ring then ${\mathcal F}$ is
Cohen-Macaulay.
\end{Proposition}
\demo Without loss of generality, we may assume $I \not= J$. As
$e_0=\lambda(R/J)$, from the identity $\lambda(I/{\mathfrak
m}I)+\lambda({\mathfrak m}I/{\mathfrak m}J) =
\lambda(I/J)+\lambda(J/{\mathfrak m}J)$ it readily follows that $-
e_0+\lambda(R/I)+\mu(I)-d = -\lambda({\mathfrak m}I/{\mathfrak
m}J)$. Thus
\[
f_0 \leq e_1 - e_0 + \lambda(R/I) + \mu(I) - d + 1 = e_1
-\lambda({\mathfrak m}I/{\mathfrak m}J)+1 \leq e_1 + 1,
\]
as desired. Therefore, the equality $f_0=e_1+1$ forces ${\mathfrak
m}I={\mathfrak m}J$ and hence ${\mathfrak m} = J \colon I$. If, in
addition, $R$ is a local Gorenstein ring, we have that $I=J \colon
{\mathfrak m}$ as well. Thus $I^2=JI$ by \cite[2.2]{CP}. Finally,
the Cohen-Macaulayness of ${\mathcal F}$ follows from
\cite[3.3]{HSa}. \QED

\medskip

We notice that equality in $($\ref{e0&e1}$)$ --- and a
fortiori in the bound established in Theorem~\ref{sally} --- does
not imply, in general, the Cohen-Macaulayness of ${\mathcal F}$
if the ring $R$ is not Gorenstein.

\begin{Example}\label{minors}
{\rm Given a field $k$, let $R$ be the power series ring
$k[\![T_1,T_2,T_3,T_4,T_5]\!]$ modulo the ideal generated by the
two by two minors of the matrix
\[
\varphi=\left[ \begin{array}{cccc} T_1 & T_2 & T_3 & T_4 \\ T_2 &
T_3 & T_4 & T_5 \end{array}\right].
\]
The ring $R$ is a local Cohen-Macaulay ring of dimension two and
type three. Let $t_i$ denote the image of $T_i$ in $R$ and
consider the ideal $I=(t_1, t_2, t_4, t_5)$. Then
\[
f_0=4=e_0=e_1+1=e_1-e_0+\lambda(R/I) + \mu(I)-d+1.
\]
However, the special fiber ring ${\mathcal F}$ of $I$ has depth
one. }
\end{Example}

\medskip

Next, we give a very general estimate for the multiplicity of
${\mathcal F}$, which is valid for any ideal $I$ of a Noetherian
local ring $(R, {\mathfrak m})$ with infinite residue field. We
observe that a characterization of the Cohen-Macaulayness of
${\mathcal F}$ in terms of the degree of ${\mathcal F}$ can also
be found in \cite[5]{Shah2}.

\begin{Proposition} \label{fiiscm}
Let $(R,{\mathfrak m})$ be a Noetherian local ring with infinite
residue field and let $I$ be an ideal with a minimal reduction
$J$. If $r=r_J(I)$ then
\[
f_0 \leq 1 + \sum_{j=1}^r \mu(I^j/JI^{j-1}),
\]
and equality holds if and only if ${\mathcal F}$ is
Cohen-Macaulay.
\end{Proposition}
\demo Set $A = {\mathcal F}(J)= k[x_1, \ldots, x_{\ell}]$, where
$\ell$ is the analytic spread of $I$, and ${\mathcal F}(I)=
\displaystyle\bigoplus_{j\geq 0}F_j$. Let $N$ be the minimal
number of generators of ${\mathcal F}(I)$ as a $A$-module. As
$f_0=\deg {\mathcal F}(I) \leq N \deg A=N$, we only need to
estimate $N$. As an $A$-module, ${\mathcal F}(I)$ is minimally
generated by
\[
\dim_k({\mathcal F}(I)/A_{+}{\mathcal F}(I)) = 1 + \sum_{j=1}^r
\dim_k(F_j/(x_1, \ldots, x_{\ell})F_{j-1})
\]
elements. But these summands are $\mu(I^j/JI^{j-1})$. The second
assertion follows from the standard criterion for
Cohen-Macaulayness as applied to graded $A$-modules. \QED

\medskip

Note that the ideal $I$ of Example~\ref{minors} has reduction
number $2$ with respect to the ideal $J=(t_1, t_2)$. Moreover,
$\mu(I/J)=2$ and $\mu(I^2/JI)=2$. Thus, the failure of the
Cohen-Macaulayness of ${\mathcal F}$ is explained by the strict
inequality $f_0=4 < 5 = 1 + \mu(I/J) + \mu(I^2/JI)$.

\medskip

Even though the equality in the bound described in
Theorem~\ref{sally} does not assure the Cohen-Macaulayness of
${\mathcal F}$, nevertheless it yields useful information about
the associated primes of ${\mathcal F}$.

\begin{Theorem} \label{cmfif}
Let $(R,\mathfrak{m})$ be a Cohen-Macaulay local ring of dimension
$d>0$ and infinite residue field. Let $I$ be an
$\mathfrak{m}$-primary ideal. If $f_0 = e_1-e_0 + \lambda(R/I)+
\mu(I)-d+1$, then ${\mathcal F}$ is unmixed.
\end{Theorem}
\demo We may assume that the reduction number of $I$ is strictly
greater than $1$. From the sequence $(\ref{tensor})$ in the proof
of Theorem~\ref{sally}, we obtain the following diagram
\[
\begin{array}{cccccccccccc}
& & & & & & & & 0 & & \\
& & & & & & & & \downarrow & & \\
& & & & & & & & K_2 & & \\
& & & & & & & & \downarrow & & \\
& & & & & & & & S_J(I)[-1] & & \\
& & & & & & & & \downarrow & & \\
0 & \!\! \rightarrow \!\! & K_1 & \!\! \longrightarrow \!\! &
{\mathcal F}(J) \oplus {\mathcal F}(J)^{n-d}[-1] & \!\!
\longrightarrow \!\! & {\mathcal F}(I) & \!\! \longrightarrow \!\!
& S_J(I)[-1] \otimes R/{\mathfrak m} & \!\!
\rightarrow \!\! & 0 \\
& & & & & & & & \downarrow & & \\
& & & & & & & & \,\,0. & &
\end{array}
\]
Now, the asserted equality $f_0=e_1-e_0+\lambda(R/I)+n-d+1$
implies that both $K_1$ and $K_2$ have multiplicity $($degree$)$
zero. But $K_1$ is a submodule of the free ${\mathcal
F}(J)$-module ${\mathcal F}(J) \oplus {\mathcal F}(J)^{n-d}[-1]$,
so it is zero. Whereas any non-zero submodule of $S_J(I)$ must
have positive multiplicity, since $S_J(I)$ has no associated
primes of height $\geq 2$ $($see \cite[2.2]{Vas}$)$. Thus $K_2$ is
zero as well. Hence we have an exact sequence
\begin{equation}\label{package}
0 \rar \mathcal{F}(J) \oplus \mathcal{F}(J)^{n-d}[-1] \lar
\mathcal{F}(I) \lar S_J(I)[-1] \rar 0,
\end{equation}
from which the desired conclusion follows. \QED

\medskip

Next, we observe that the previous theorem gives us a relation
between the depths of ${\mathcal F}$ and $S_J(I)$, and a fortiori
a relation between the depths of ${\mathcal F}$ and ${\mathcal
G}$. In particular, good depth conditions on ${\mathcal G}$ and
the equality in the bound established in Theorem~\ref{sally} force
the Cohen-Macaulayness of ${\mathcal F}$.

\begin{Corollary}\label{depthd-1}
Let $(R,\mathfrak{m})$ be a Cohen-Macaulay local ring of dimension
$d>0$ and infinite residue field. Let $I$ be an
$\mathfrak{m}$-primary ideal. If $f_0 = e_1-e_0 + \lambda(R/I)+
\mu(I)-d+1$, then
\[
{\rm depth}\, {\mathcal F} \geq \min \{ {\rm depth}\, {\mathcal
G}+1, d \}.
\]
If, in addition, ${\rm depth}\, {\mathcal G} \geq d-1$ then
${\mathcal F}$ is Cohen-Macaulay.
\end{Corollary}
\demo The inequality follows from a depth count on
$(\ref{package})$ and the equality ${\rm depth}\, S_J(I) = \min \{
{\rm depth}\, {\mathcal G}+1, d \}$, which follows from
\cite[1.2.11]{M-thesis} and \cite[2.1]{Vaz}. \QED

\bigskip

The associated primes of $\mathcal{F}$ play a role in the integral
closedness of the product $\mathfrak{m}I$. To see how this comes
about, we recall an observation of \cite[1.5]{HuHu}, later refined
in \cite[4.1]{CGPU}. Then, in Corollary~\ref{integrality}, we give
an application in its spirit.

\begin{Proposition}
Let $(R, \mathfrak{m})$ be a normal local domain of dimension
$d>0$ with infinite residue field and suppose $I$ is a normal
ideal of analytic spread $d$. If $\mathcal{F}$ is unmixed then
$\overline{\mathfrak{m}I^m} = {\mathfrak m}I^m$ for all $m$.
\end{Proposition}

\begin{Corollary}\label{integrality}
Let $(R,\mathfrak{m})$ be a normal Cohen-Macaulay local domain of
dimension $d>0$ and infinite residue field. Let $I$ be an
$\mathfrak{m}$-primary normal ideal. If $f_0 = e_1-e_0 +
\lambda(R/I)+ \mu(I)-d+1$, then $\overline{\mathfrak{m}I^m} =
{\mathfrak m}I^m$ for all $m$.
\end{Corollary}

\medskip

In the same fashion as Theorem~\ref{sally} and Theorem~\ref{cmfif},
one can show a bound on the first two Hilbert coefficients of $I$,
which possibly yields information on the associated graded ring
${\mathcal G}$ of $I$. The following inequality is of interest in
the case of exotic local Cohen-Macaulay rings of small dimension,
that are not regular local rings. The bound was proved by Elias
and Valla in the case of the maximal ideal \cite[Section 2]{EliasValla}.

\begin{Remark}
{\rm Let $(R, {\mathfrak m})$ be a local Cohen-Macaulay ring of
dimension $d>0$ and infinite residue field. Let $I$ be an
${\mathfrak m}$-primary ideal. Then the coefficients
$e_0$ and $e_1$ of the Hilbert polynomial of $I$ satisfy
\[
2e_0 - e_1 \leq \lambda(R/I) (\mu(I) - d + 2).
\]
Moreover, if equality holds one has that the associated graded ring
${\mathcal G}$ is unmixed; if, in addition, the Sally module $S_J(I)$ of
$I$ with respect to a minimal reduction $J$ is Cohen-Macaulay
then ${\mathcal G}$ is Cohen-Macaulay as well.}
\end{Remark}

\medskip

Another interesting example, where we can argue the Cohen-Macaulayness
of ${\mathcal F}$ and hence the integral
closedness of the ideals ${\mathfrak m}I^m$ for all $m$, is due to
Huckaba and Huneke \cite[3.11]{HH}. It is the first explicit example of
the failure in dimension three of a two-dimensional version $($due
to Sancho de Salas$)$ of a vanishing theorem of Grauert and
Riemenschneider. This example complements nicely a family of such
examples over ${\mathbb C}$ constructed by Cutkosky \cite[Section III]{C}.

\begin{Example}\label{ex-HuHu}{\rm
Let $k$ be a field of characteristic $\not= 3$ and set $R =
k[\![x,y,z]\!]$, where $x,y,z$ are indeterminates. Let $P = (x^4,
x(y^3+z^3), y(y^3+z^3), z(y^3+z^3))$ and set $I = P + {\mathfrak
m}^5$, where ${\mathfrak m}$ is the maximal ideal of $R$. The
ideal $I$ is a $16$ generated, normal, ${\mathfrak m}$-primary
ideal with ${\rm depth}\,{\mathcal G}(I^m) = {\rm dim}\,R-1(=2)$
for all $m \geq 1$.  Moreover, we have that $I^4=JI^3$, where $J=
(x^4, z(y^3+z^3), y(y^3+z^3)+z^5)$. The Cohen-Macaulayness of
${\mathcal F}$ follows from Proposition~\ref{fiiscm} since we
have that $f_0=16$, $\mu(I/J)=13$ and $\mu(I^2/JI)=
\mu(I^3/JI^2)=1$. On the other hand, as far as the bound in
Theorem~\ref{sally} is concerned, we have that $f_0 = 16 < 17 =
e_1 - e_0 + \lambda(R/I) + \mu(I) -d + 1$. In fact one has that
$e_0 = 76$, $e_1 = 48$ and $\lambda(R/I)=31$. }
\end{Example}

\subsection{Lower bounds}

Let $(R, {\mathfrak m})$ be a Noetherian local ring and let $I$ be
an ${\mathfrak m}$-primary ideal. Let $R/I=A_0 \supset A_1 \supset
\cdots \supset A_s=0$, with $\lambda(A_i/A_{i+1})=1$ for $0 \leq i
\leq s-1$ and $s=\lambda(R/I)$, be a composition series of $R/I$. Thus
we can find $u_i \in R$ such that $A_i=(A_{i+1}, u_i)$ and ${\mathfrak
m}\,u_i\in A_{i+1}$. On the other hand, this composition series
induces a filtration ${\mathcal G}=A_0 {\mathcal G} \supset A_1
{\mathcal G} \supset \cdots \supset A_s {\mathcal G}=0$ on the
associated graded ring ${\mathcal G}$ of $I$, whose factors
$A_i{\mathcal G}/A_{i+1}{\mathcal G} \cong u_i{\mathcal G}$ are
${\mathcal F}(I)$-modules. A bookkeeping in the family of
epimorphisms ${\mathcal F} \longrightarrow u_i{\mathcal G}
\rightarrow 0$ yields the bound $e_0 \leq \lambda(R/I) \, f_0$,
or, equivalently,
\[
\frac{e_0}{\lambda(R/I)} \leq f_0.
\]
It has been proved in \cite[2.6]{HKU} that if the associated
primes of ${\mathcal F}$ have the same dimension then equality
in the above formula is equivalent to $I^m/I^{m+1}$ being $R/I$-free
for all $m \geq 1$. In this situation, $R$ is said to be normally
flat along $I$. We also observe that an earlier version of this
latter characterization can be found in \cite[8]{Shah2}.

\section{Special bounds}

In \cite[2.2]{S0}, Sally proved that the reduction number of the
maximal ideal of a local Cohen-Macaulay ring of dimension $d$ is
bounded by $d!\, e_0-1$. In cite \cite[6.12 and 6.16]{deg}, Vasconcelos
improved this bound and for any ${\mathfrak m}$-primary ideal he
showed that $r \leq d(e_0-2)+1$. In this Section we give bounds on
the reduction number of an ${\mathfrak m}$-primary ideal via the
estimates on the multiplicity of the special fiber ring ${\mathcal
F}$ obtained earlier.

\begin{Corollary}\label{firstcor}
Let $(R,\mathfrak{m})$ be a Cohen-Macaulay local ring of dimension
$d>0$. Let $I$ be an $\mathfrak{m}$-primary ideal with reduction
number $r$. If the residue field $R/{\mathfrak m}$ has
characteristic $0$ and $f_0 = e_1 - e_0 + \lambda(R/I)+
\mu(I)-d+1$, then
\[
r \leq e_1 - e_0 +\lambda(R/I) + \mu(I) -d.
\]
\end{Corollary}
\demo From Theorem~\ref{cmfif} we conclude that ${\mathcal F}$ is
unmixed. Hence the reduction number of $I$ is at most $f_0-1$, by
a result of Vasconcelos \cite[7 and 9]{compositio} $($see also
\cite[2.2$(${\it a}$)$]{reds}$)$. \QED

\medskip

We note that the above bound on $r$ can actually be sharpened if
the special fiber ring ${\mathcal F}$ is Cohen-Macaulay: See
Corollary~\ref{coro-rossi} below.

\begin{Remark}\label{CMfiber}
{\rm Let $R$ be a local Cohen-Macaulay ring of dimension $d>0$ and
infinite residue field. Let $I$ be an ${\mathfrak m}$-primary
ideal of reduction number $r$. If the special fiber ring
${\mathcal F}$ of $I$ is Cohen-Macaulay then
\[
r \leq f_0 - \mu(I) + d,
\]
where $f_0$ is the multiplicity of ${\mathcal F}$.}
\end{Remark}
\demo We seek to estimate the degree $r$ of the last basis element
of the vector space ${\mathcal V} = {\mathcal F}/({\bf
z}){\mathcal F}$, where ${\bf z}=z_1, \ldots, z_d$ is a
superficial sequence. As ${\mathcal V}$ has one basis element in
degree zero, $\mu(I)-d$ basis elements in degree one and at least
one basis element in each degree,
we get $1 + \mu(I) - d + r-1 \leq {\rm dim}\, {\mathcal V} = f_0$.
Hence our claim. \QED

\medskip

The assumption on the Cohen-Macaulayness of ${\mathcal F}$ is
necessary as the next example shows.

\begin{Example}
{\rm Let $R=k[\![t^6,t^{11},t^{15},t^{31}]\!]$, where $k$ is an
infinite field, and let $I=(t^6,t^{11},t^{31})$. The ideal $I$ has
height 1, analytic spread 1 and reduction number 2. Furthermore
${\mathcal G}(I)$ is Cohen-Macaulay, but ${\mathcal F}$ is not
Cohen-Macaulay. One has that \ $e_0=6, \ e_1 = 5$ \ and \ $f_0=3$.
Thus $r=2 > 1 = f_0-\mu(I)+d$. On the other hand one has
$r=2=e_1-e_0+\lambda(R/I)+1$.}
\end{Example}

\medskip

Next, we point out a consequence of Remark~\ref{CMfiber} and
Theorem~\ref{sally} in the special case of ideals with reduction
number one. It easily recovers an earlier result of Shah
\cite[7$({\it b})$]{Shah2}.

\begin{Corollary}\label{r=1}
Let $R$ be a local Cohen-Macaulay ring of dimension $d>0$ and
infinite residue field. Let $I$ be an ${\mathfrak m}$-primary
ideal of reduction number one. Then the multiplicity $f_0$ of the
special fiber ring ${\mathcal F}$ of $I$ is given by
$f_0=\mu(I)-d+1$.
\end{Corollary}
\demo As $I$ has reduction number one it follows from
\cite[2.1]{H} and \cite[3.2 and 3.3]{O}
that $e_1-e_0+\lambda(R/I)=0$. Thus, Theorem~\ref{sally} gives
that $f_0 \leq \mu(I)-d+1$. On the other hand, the special fiber
ring of $I$ is Cohen-Macaulay by \cite[3.3]{HSa}, so that the
reverse inequality follows from Abhyankar's bound \cite[1]{Abh}.
\QED

\medskip

The following result has been proved without any assumption on the
special fiber ring by Rossi for local Cohen-Macaulay rings of
dimension at most two \cite[1.5]{R}. It also holds whenever ${\rm
depth}\, {\mathcal G} \geq {\rm dim}\, R-1$.

\begin{Corollary}\label{coro-rossi}
Let $R$ be a local Cohen-Macaulay ring of dimension $d>0$ and
infinite residue field. Let $I$ be an ${\mathfrak m}$-primary
ideal of reduction number $r$. If the special fiber ring
${\mathcal F}$ of $I$ is Cohen-Macaulay then
\[
r \leq e_1-e_0+\lambda(R/I)+1.
\]
\end{Corollary}
\demo By combining Remark~\ref{CMfiber} and Theorem~\ref{sally}
one has that
\begin{eqnarray*}
r & \leq & f_0 - \mu(I) + d \hspace{6.3cm} (\mbox{by Remark~\ref{CMfiber}}) \\
  & \leq & (e_1 - e_0 + \lambda(R/I) + \mu(I) - d + 1) - \mu(I) + d
\hspace{1cm} (\mbox{by Theorem~\ref{sally}}) \\
  &   =  & e_1 - e_0 + \lambda(R/I) + 1,
\end{eqnarray*}
as claimed. \QED

\medskip

If $R$ is a $2$-dimensional regular local ring, we obtain an
estimate on the reduction number $r$ and on the multiplicity $f_0$
of an ideal $I$ which does not explicitly involve the Hilbert
coefficients of the ideal $I$.

\begin{Corollary}
Let $(R, {\mathfrak m})$ be a $2$-dimensional regular local ring
of infinite residue field. Let $I$ be an ${\mathfrak m}$-primary
ideal of reduction number $r$. Then
\[
r \leq \lambda(\overline{I}/I)+1  \qquad {\it and} \qquad f_0 \leq
\lambda(\overline{I}/I)+\mu(I)-1,
\]
where $\overline{I}$ denotes the integral closure of $I$.
\end{Corollary}
\demo We have that $e_1(I) \leq \overline{e}_1(I)$ by
\cite{PUVa}, where $\overline{e}_1(I)$ is the corresponding Hilbert
coefficient arising from the filtration given by the integral closure of the
powers of $I$ rather than the $I$-adic filtration.
On the other hand, using the fact $\overline{I}$ is
a normal ideal of reduction number one \cite[5.5]{LT}, one has that
$\overline{e}_1(I) = e_1(\overline{I}) = \lambda(\overline{I}/J)$
for any   minimal reduction $J$ of $I$. As $e_0=\lambda(R/J)$, one
finally concludes from Rossi's bound \cite[1.5]{R} that
\[
r \leq e_1-e_0+\lambda(R/I)+1 \leq
\lambda(\overline{I}/J)-e_0+\lambda(R/I)+1=\lambda(\overline{I}/I)+1,
\]
as claimed. \QED

\medskip

\begin{Remark}{\rm
Observe that whenever $\overline{I} \not= {\mathfrak m}$ one has
that the above bounds translate into
\[
r \leq e_0-2 \qquad  {\rm and} \qquad f_0 \leq
\lambda(R/I)+\mu(I)-3.
\]
The first one compares favorably with Vasconcelos' earlier bound,
which says $r \leq 2e_0 -3$.}
\end{Remark}

\section{Partially identical Hilbert polynomials}

To compare $I$ to its integral closure $\overline{I}$ it is useful
to track their numerical measures as expressed by appropriate
Hilbert functions. An approach was given in \cite[1]{Shah} to
construct a `canonical' sequence of ideals $I \subset I_0 \subset
I_1 \subset \cdots \subset I_{d-1}\subset I_d = \overline{I}$,
with Hilbert polynomials partially identical to the one of $I$.
More recently, Ciuperca \cite[2.5]{Ciu} has instead shown a
connection with the $S_2$-ification of the Rees algebra of $I$.

We give below another derivation of Shah's result as it will be
relevant for an improvement of our earlier bounds on $f_0$.

\begin{Theorem} \label{s2ofideal}
Let $(R, {\mathfrak m})$ be an analytically equidimensional local
domain of dimension $d$ and let $I$ be an ${\mathfrak m}$-primary
ideal.
There exists
a unique, largest ideal $I\subset I_j \subseteq \overline{I}$, with
$0 \leq j \leq d$ and the property that the corresponding Hilbert
polynomial has coefficients satisfying
\[
e_i(I_j) = e_i(I), \quad i=0,\ldots, d-j.
\]
\end{Theorem}
\demo Let $A$ be the extended Rees algebra $R[It,t^{-1}]$ of $I$.
For each integer $1\leq j \leq d+1$, let $B^{(j)}$ be the subring
of $C=R[t,t^{-1}] \cap \overline{R[It,t^{-1}]}$ defined by
$B^{(j)}= \{  h \in C \ | \ {\rm ht}\, (A :_A h) \geq j \}$. This
gives a filtration of subalgebras $A=B^{(d+2)}\subset B^{(d+1)}
\subset \cdots \subset B^{(2)} \subset B^{(1)}$. By considering
only the component $A^{(j)}$ of $B^{(j)}$ in non-negative degrees
$($since they all have the same components in negative degrees$)$
we have
\[
\begin{array}{ccl}
A^{(1)} &=& \text{integral closure of $R[It]$ in
$R[t]$} \\
A^{(2)} &=& \text{$S_2$-ification of $R[It]$}\\
A^{(d+1)}  &=& \text{Ratliff-Rush closure of $R[It]$}\\
A^{(d+2)}  & =& \text{Rees algebra $R[It]$ of $I$}.
\end{array}
\]

We are now ready to relate these algebras to the ideals $I_j$ and
to prove the assertion about their Hilbert polynomials. Set $I_j$
to be the ideal defined by the component of degree $1$ of
$A^{(d+1-j)}$. We have exact sequences of finitely generated
modules over $R[It,t^{-1}]$
\begin{eqnarray*} \label{fv1}
0 \rar R[It,t^{-1}]\lar R[I_jt,t^{-1}]\lar M_j\rar 0,
\end{eqnarray*}
where $M_j$ is a module of dimension at most $d+1-(d+1-j)=j$.

Tensoring this sequence by $R[It,t^{-1}]/(t^{-1})$, we obtain two
exact sequences of graded modules:
\begin{eqnarray}\label{fv2}
0 \rar K_j[+1]\lar M_j[+1] \stackrel{t^{-1}}{\lar} M_j\lar   C_j
\rar 0,
\end{eqnarray}
\begin{eqnarray}\label{fv3}
0 \rar K_j \lar {\mathcal G}(I) \lar {\mathcal G}(I_j) \lar C_j
\rar 0,
\end{eqnarray}
Note that $($\ref{fv2}$)$ was induced by multiplication by
$t^{-1}$, an endomorphism which is nilpotent on $M_j$; it follows
that (see \cite[12.1]{Eisenbudbook})  the Krull dimension of all
the modules in this sequence is the same.  One can further assert
that the multiplicities of $K_j$ and of $C_j$ match. As a
consequence, adding in (\ref{fv3}) the Hilbert polynomials, we
obtain that the Hilbert polynomials of ${\mathcal G}(I)$ and of
${\mathcal G}(I_j)$ have matching coefficients not only down to
degree $j$ $($which is guaranteed by the codimension of $M_j)$ but
one position further down. \QED

\medskip

One can actually do slightly better in the estimate of $f_0$
established in Theorem~\ref{sally}.

\begin{Corollary}\label{rr-closure}
Let $R$ be a local Cohen-Macaulay ring of dimension $d$ and
infinite residue field. Let $I$ be an ${\mathfrak m}$-primary
ideal. Then the multiplicity $f_0$ of the special fiber ring
${\mathcal F}$ of $I$ satisfies
\begin{equation}\label{expmultfi}
f_0 \leq e_1 - e_0 + \lambda(R/\widetilde{I}) + \mu(\widetilde{I})
- d + 1,
\end{equation}
where $\widetilde{I}$ is the Ratliff-Rush closure of $I$.
\end{Corollary}
\demo The special fiber rings ${\mathcal F}(\widetilde{I})$ and
${\mathcal F}(I)$ have the same multiplicity, since these two
algebras differ at most in a finite number of components. Now
observe that the relations
\[
\lambda(R/\widetilde{I})=\lambda(R/I)-\lambda(\widetilde{I}/I)
\qquad {\rm and} \qquad \mu(\widetilde{I}) \leq \mu(I)+
\lambda(\widetilde{I}/I)
\]
implies that $\lambda(R/\widetilde{I}) + \mu(\widetilde{I}) \leq
\lambda(R/I) + \mu(I)$. Finally, since the Hilbert coefficients
are the same for both ideals $\widetilde{I}$ and $I$, we conclude
that the formula
\[
f_0 \leq e_1 - e_0 + \lambda(R/\widetilde{I}) + \mu(\widetilde{I})
- d + 1
\]
sharpens Theorem~\ref{sally}. \QED

\medskip

An additional improvement occurs from the following consideration.
Let $\check{I}$ denote the largest ideal containing $I$ such that
in the embedding of Rees algebras
\begin{eqnarray}\label{eq8}
0\rar {\mathcal R}(I) \lar {\mathcal R}(\check{I}) \lar C \rar 0,
\end{eqnarray}
${\rm dim}\, C \leq d-1$. Thus, $\check{I}$ is the ideal $I_{d-1}$
in the terminology established in the proof of
Theorem~\ref{s2ofideal}, that is the degree one component of the
$S_2$-ification of the Rees algebra of $I$.

\begin{Corollary}\label{e1-closure}
Let $R$ be a local Cohen-Macaulay ring of dimension $d>0$ and
infinite residue field. Let $I$ be an ${\mathfrak m}$-primary
ideal. Then the multiplicity $f_0$ of the special fiber ring
${\mathcal F}$ of $I$ satisfies
\[
f_0 \leq e_1 - e_0 + \lambda(R/\check{I}) + \mu(\check{I}) - d +1.
\]
\end{Corollary}
\demo According to Theorem~\ref{s2ofideal} $\check{I}$ is also the
largest ideals containing $I$ with the same values for the Hilbert
coefficients $e_0$ and $e_1$. From equation $(\ref{eq8})$ it
follows that there is an induced sequence
\[
0 \rar D \lar {\mathcal F}(I) \lar {\mathcal F}(\check{I})\lar
C/{\mathfrak m}C \rar 0
\]
of finitely generated ${\mathcal F}(I)$-modules such that the
dimensions of $C/{\mathfrak m}C$ and $D$ are at most $d-1$.
Indeed, in the case of $D$ one has the map ${\rm Tor}_1(C,
R/{\mathfrak m}) \lar D \rar 0$, which justifies our assertion
about the dimension. Thus, the special fiber rings ${\mathcal
F}(\check{I})$ and ${\mathcal F}(I)$ have the same multiplicity.
Therefore if we write formula (\ref{expmultfi}) for $\check{I}$
instead of $\widetilde{I}$,
\[
f_0 \leq e_1 - e_0 + \lambda(R/\check{I}) + \mu(\check{I})
- d + 1
\]
we can only achieve gains, as the function
$\lambda(R/\cdot)+\mu(\cdot)$ is monotone non-increasing. \QED

\medskip

\begin{Corollary}\label{mu-1}
Let $R=k[x,y]$ with $k$ a field of characteristic zero and $x, y$
variables, and let ${\mathfrak m}$ denote the maximal homogeneous
ideal of $R$. Let $I$ be an ${\mathfrak m}$-primary monomial ideal
with a monomial $2$-generated reduction. Then the multiplicity
$f_0$ of the special fiber ring ${\mathcal F}$ of $I$ is given by
$f_0 = \mu(\check{I})-1$, where $\check{I}$ is the degree one
component of the $S_2$-ification of the Rees algebra of $I$.
\end{Corollary}
\demo For any 2-generated reduction $J$ of $I$, it follows from
\cite{PUVi} that $\check{I}=J \colon (J^r\colon I^r)$, where $r$
is the reduction number of $I$ with respect to $J$. Moreover,
$\check{I}$ has reduction number one. Our assertion then follows
from Corollary~\ref{r=1}. \QED

\medskip

The assertion of Corollary~\ref{mu-1} seems to hold for all
${\mathfrak m}$-primary ideals of a two-dimensional regular local
ring. The following example is taken from \cite[3.3]{Ciu}, where
it is erroneously stated that ${\mathcal R}(\check{I})$ is not
Cohen-Macaulay.

\begin{Example}{\rm
Let $R=k[x,y]_{(x,y)}$ with $k$ a field of characteristic zero.
The ideal $I=(x^8,x^3y^2,x^2y^4,y^8)$ has $S_2$-ification
$\check{I}=(x^8, x^3y^2, x^2y^4, xy^6, y^8)$. We have that $f_0=4,
\,e_0=40$ and $e_1=12$ whereas $\lambda(R/I)=30$ and
$\lambda(R/\check{I})=28$. Thus in one case we have that $f_0=4 <
5=e_1-e_0+\lambda(R/I)+\mu(I)-d+1$ while in the other we have that
$f_0=4=e_1-e_0+\lambda(R/\check{I})+\mu(\check{I})-d+1$. By
choosing $J=(x^8+y^8+x^2y^4, x^3y^2)$ we have that
$\check{I}^2=J\check{I}$ so that $f_0 = \mu(\check{I})-1$.}
\end{Example}

\medskip

Motivated by Theorem~\ref{s2ofideal} and its proof, we can also
show that the same sequence of ideals $I_j$, with $0 \leq j \leq
d-1$, has special fiber rings ${\mathcal F}(I_j)$ with partially
identical Hilbert polynomials to the one of ${\mathcal F}(I)$.
Consistently with our previously established terminology, we
denote with $f_0=f_0(I), f_1=f_1(I), \ldots, f_{d-1}=f_{d-1}(I)$
the coefficients of the Hilbert polynomial associated with the
Hilbert function of the special fiber ring ${\mathcal F}(I)$.

\begin{Theorem}
Let $(R, {\mathfrak m})$ be an analytically equidimensional local
domain of dimension $d$ and let $I$ be an ${\mathfrak m}$-primary
ideal. There exists a unique, largest ideal $I\subset I_j
\subseteq \check{I}$, with $0 \leq j \leq d-1$ and the property
that the corresponding Hilbert polynomial of the special fiber
rings has coefficients satisfying
\[
f_i(I_j) = f_i(I), \quad i=0,\ldots, d-j-1.
\]
\end{Theorem}

In particular, the above result $($see also the proof of
Corollary~\ref{e1-closure}$)$ says that the degree one component
$\check{I}$ of the $S_2$-ification of the Rees algebra ${\mathcal
R}$ of $I$ can also be characterized as the {\it largest} ideal
containing $I$ with the same multiplicity $f_0$. This statement is
analogous to the characterization of the integral closure
$\overline{I}$ of $I$ as the largest ideal containing $I$ with the
same multiplicity $e_0$. However, our interest in $\check{I}$ has
also been motivated by a recent connection, made in \cite{PUVi},
with the computation of the {\it core} of $I$ $($we recall that
${\rm core}(I)$ is the intersection of all the reductions of $I)$.
The core of an ideal has lately been under much scrutiny: One of
the classical motivations to study the core comes from the
Brian\c{c}on-Skoda theorem, but more recently Hyry and Smith have
shown that Kawamata's well-known conjecture on the existence of
sections of line bundles is equivalent to a statement about the
core of certain ideals in section rings. In a more general setting
that ours, Polini, Ulrich and Vitulli have shown that ${\rm
core}(I) = {\rm core}(\check{I})$. In view of our result, their
statement can also be formulated by saying that two ${\mathfrak
m}$-primary ideals $I \subset K$ with the same multiplicity $f_0$
also satisfy ${\rm core}(I)={\rm core}(K)$.


\bigskip

\begin{thebibliography}{99}

\bibitem[{\bf 1}]{Abh}{S.S. Abhyankar, Local rings of high embedding
dimension, Amer. J. Math. {\bf 89} (1967), 1073-1077.}

\bibitem[{\bf 2}]{Ciu}{C. Ciuperca, First coefficient ideals and the
$S_2$-ification of a Rees algebra, J. Algebra {\bf 242} (2001),
782-794.}

\bibitem[{\bf 3}]{CP}{A. Corso and C. Polini, Links of prime
ideals and their Rees algebras, J. Algebra {\bf 178} (1995),
224-238.}


\bibitem[{\bf 4}]{CGPU}{A. Corso, L. Ghezzi, C. Polini and B. Ulrich,
Cohen-Macaulayness of special fiber rings, Comm. Algebra
$($Kleiman's volume$)$ {\bf 31} (2003), 3713-3734.}

\bibitem[{\bf 5}]{C}{S.D. Cutkosky, A new characterization of rational
surfaces singularities, Invent. Math. {\bf 102} (1990) 157-177.}


\bibitem[{\bf 6}]{Eisenbudbook}{D. Eisenbud, Commutative algebra with a view
toward algebraic geometry, Springer-Verlag, New York, 1995.}


\bibitem[{\bf 7}]{EliasValla}{J. Elias and G. Valla, Rigid Hilbert
functions, J. Pure and Appl. Algebra {\bf 71} (1991), 19-41.}

\bibitem[{\bf 8}]{HKU}{W. Heinzer, M.-K. Kim and B. Ulrich, The Gorenstein and
complete intersection properties of associated graded rings, to
appear in J. Pure and Appl. Algebra $($Vasconcelos' volume$)$.}


\bibitem[{\bf 9}]{HuHu}{R. H\"{u}bl and C. Huneke, Fiber cones and the
integral closure of ideals, Collect. Math. {\bf 52} (2001),
85-100.}

\bibitem[{\bf 10}]{HH}{S. Huckaba and C. Huneke, Normal ideals in
regular rings, J. reine angew. Math. {\bf 510} (1999), 63-82.}

\bibitem[{\bf 11}]{H}{C. Huneke, Hilbert functions and symbolic
powers, Michigan Math. J. {\bf 34} (1987), 293-318.}


\bibitem[{\bf 12}]{HSa}{C. Huneke and J.D. Sally, Birational extensions in
dimension two and integrally closed ideals, J. Algebra {\bf 115}
(1988), 481-500.}




\bibitem[{\bf 13}]{LT}{J. Lipman and B. Teissier,  Pseudo-rational local rings
and a theorem of Brian\c{c}on-Skoda, Michigan Math. J. {\bf 28}
(1981), 97-116.}

\bibitem[{\bf 14}]{Mazur}{B. Mazur, Deformations of Galois
representations and Hecke algebras, Harvard course notes, available by
request from the author, 1994.}


\bibitem[{\bf 15}]{NR}{D.G. Northcott and D. and Rees, Reductions of
ideals in local rings, Proc. Camb. Phil. Soc. {\bf 50} (1954), 145-158.}

\bibitem[{\bf 16}]{O}{A. Ooishi, $\Delta$-genera and sectional genera of
commutative rings, Hiroshima Math. J. {\bf 17} (1987), 361-372.}


\bibitem[{\bf 17}]{PUVa}{C. Polini, B. Ulrich, W.V.
Vasconcelos and R. Villarreal, Normalization of ideals, in
preparation.}

\bibitem[{\bf 18}]{PUVi}{C. Polini, B. Ulrich and M. Vitulli,
Core of monomial ideals, in preparation.}

\bibitem[{\bf 19}]{RR}{L.J. Ratliff and D.E. Rush,
Two notes on reductions of ideals, Indiana Univ. Math. J. {\bf 27}
(1978), 929-934.}

\bibitem[{\bf 20}]{R}{M.E. Rossi, A bound on the reduction number of a
primary ideal, Proc. Amer. Math. Soc. {\bf 128} (2000),
1325-1332.}

\bibitem[{\bf 21}]{S0}{J.D. Sally, Bounds for number of generators
for Cohen-Macaulay ideals, Pacific J. Math. {\bf 63} (1976),
517-520.}

\bibitem[{\bf 22}]{Shah}{K. Shah, Coefficient ideals, Trans. Amer. Math. Soc.
{\bf 327} (1991), 373-384.}

\bibitem[{\bf 23}]{Shah2}{K. Shah, On the Cohen-Macaulayness of the
fiber cone of an ideal, J. Algebra {\bf 143} (1991), 156-172.}

\bibitem[{\bf 24}]{Vas}{W.V. Vasconcelos, Hilbert functions, analytic
spread and Koszul homology, Contemporary Mathematics {\bf 159}
(1994), 410-422.}


\bibitem[{\bf 25}]{compositio}{W.V. Vasconcelos, The redution number of an
algebra, Compositio Math. {\bf 104} (1996), 189-197.}

\bibitem[{\bf 26}]{deg}{W.V. Vasconcelos, Cohomological degrees of
graded modules, in {\it Six Lectures on Commutative Algebra}, J. Elias, J.M.
Giral, R.M. Mir\'o-Roig and S. Zarzuela (eds.), Progress in
Mathematics {\bf 166}, Birkh\"auser, Basel, 1998, 345-392.}

\bibitem[{\bf 27}]{reds}{W. V. Vasconcelos, The reduction numbers of an ideal,
J. Algebra {\bf 216} (1999), 652-664.}

\bibitem[{\bf 28}]{mrn}{W.V. Vasconcelos, Multiplicities and reduction
numbers, to appear in Compositio Math.}

\bibitem[{\bf 29}]{Vaz}{M.T.R. Vaz Pinto, Structure of Sally modules and
Hilbert functions, Ph.D. thesis, 1995.}

\bibitem[{\bf 30}]{M-thesis}{M.T.R. Vaz Pinto, Hilbert functions and Sally
modules, J. Algebra {\bf 192} (1997), 504-523.}


\bibitem[{\bf 31}]{W}{A. Wiles, Modular elliptic curves and Fermat's
last theorem, Ann. of Math. {\bf 141} (1995), 443-551.}

\end{thebibliography}
\end{document}